\documentclass{amsart}
\usepackage{amstext,amsfonts,amsmath,amsthm,amssymb,amscd}

\begin{document}

\title{Plane Jacobian conjecture for rational polynomials}
\author{Nguyen Van Chau}
\address{Institute of Mathematics, 18 Hoang Quoc Viet, 10307 Hanoi,Vietnam}
\email{nvchau@math.ac.vn}

\date{}
\maketitle

\rightline{\small \it In memory of Professor Carlos Gutierrez
\hskip 0.2 cm}

\begin{abstract}
A non-zero constant Jacobian polynomial maps $F=(P,Q)$ of
$\mathbb{C}^2$ is invertible if  $P$ and $Q$ are rational
polynomials.

{\it Keywords and Phrases:} Jacobian conjecture, Rational
polynomial.

{\it 2000 Mathematical Subject Classification:} 14R15, 14H20.

\end{abstract}

% Private commands
\newtheorem{theorem}{Theorem}
\newtheorem{lemma}{Lemma}
\newtheorem{proposition}{Proposition}

\vskip 0.3cm \noindent {\bf 1.} We shall call a polynomial map
$F=(P,Q):\mathbb{C}^2 \longrightarrow \mathbb{C}^2$ is  {\it
Keller map} if $F$ satisfies the Jacobian condition
$J(P,Q):=P_xQ_y-P_yQ_x\equiv c\neq 0$. The Jacobian conjecture,
posed first in 1939 by Ott-Heinrich Keller \cite{Keller} and still
opened, asserts that {\it every Keller map is invertible}. We
refer the readers to the nice surveys \cite{Bass} and
\cite{Essen-book} for the history, the recent developments and the
related topics of this mysterious problem.

One of most simple topological cases of Keller's question, for
which we may hope to have a complete solution by using the rich
knowledge on plane algebraic curves, is the case when one or both
of $P$ and $Q$ are {\it rational polynomials}, i.e. polynomials
with generic fibre diffeomorphic to the sphere with finite number
of punctures. Since 1978 Razar had found the following.
\begin{theorem}[Razar's Theorem, \cite{Razar}] A Keller map $F=(P,Q)$ is invertible if
$P$ is a rational polynomial with all irreducible fibres.
 \end{theorem} In attempting to understand the
nature of the plane Jacobian conjecture, Razar 's theorem had been
reproved by Heitmann \cite{Heitmann}, L\^e and Weber \cite{LeWe2},
Friedland \cite{Friedland}, Nemethi and Sigray
\cite{NemethiSigray} in several different algebraic and
algebro-geometric approaches. In fact, Vistoli \cite{Vistoli} and
Neumann and Norbury \cite{NeumannNorbudy} observed that every
rational polynomial with all irreducible fibres is equivalent to
the projection $(x,y)\mapsto x$ up to algebraic coordinates.
Recently, L\^e (\cite{Le-hanoi}, \cite{Le-kyoto}) proved  that a
Keller map $F=(P,Q)$ is invertible if $P$ is a rational polynomial
and, in addition, $P$ is simple, i.e. in regular extensions
$p:X\longrightarrow \mathbb{P}^1$ of $P$ over a compactification
$X$ of $\mathbb{C}^2$ the restriction of $p$ to each irreducible
component of the divisor at infinity $\mathcal{D}:=X\setminus
\mathbb{C}^2$ have degree $0$ or $1$. As shown in
\cite{Chau2007b}, L\^e's result is still true without the
condition that $P$ is rational.

In this paper we would like to note that  the plane Jacobian
conjecture is true for the case when both of $P$ and $Q$ are
rational.

\begin{theorem}[Main Theorem]\label{Main}
Suppose $F=(P,Q)$ is a Keller map. If $P$ and $Q$ are rational
polynomials, then $F$ is invertible.
\end{theorem}

By a {\it polar branch} we mean an irreducible branch curve at
infinity along which $F$ tends to infinity. An obvious simple fact
is that under the Jacobian condition any irreducible component of
any fiber of $P$ must contains some polar branches. Otherwise, the
restriction of $F$ to such an exceptional component must be
constant mapping that is impossible. In order to prove Theorem
\ref{Main} we will try to show that each of fibres of $P$ has only
one polar branch. This ensures that all fibres of $P$ are
irreducible. Then,  invertibility of $F$ follows from Razar's
theorem.

\vskip0.3cm \noindent{\bf 2. } Our proof  is based on the
following facts on Keller maps.

Following \cite{Jelonek}, {\it the non-proper value set } $A_f$ of
a polynomial map $f$ of $\mathbb{C}^2$ is the set of all values
$a\in \mathbb{C}^2$ such that $ f(b_i)\rightarrow a$ for a
sequence $b_i\rightarrow \infty$. This set $A_f$ is either empty
or an algebraic curve in $\mathbb{C}^2$ each of whose irreducible
components is the image of a non-constant polynomial map from
$\mathbb{C}$ into $\mathbb{C}^2$. If $f$ is a Keller map, the
restriction
$$f: \mathbb{C}^2\setminus f^{-1}(A_f)\longrightarrow \mathbb{C}^2\setminus A_f$$
gives a unbranched covering.

\begin{theorem}[ see Theorem 4.4 in \cite{Chau99}]\label{AF} Suppose $F=(P,Q)$ is a Keller map. The non-proper value set $A_F$, is not empty, is
composed of the images of some polynomial maps $ t\mapsto
(\alpha(t),\beta(t))$, $\alpha, \beta\in \mathbb{C}[t]$,
satisfying $${\deg \alpha\over\deg \beta}={\deg P\over\deg
Q}.\eqno(1)$$ In particular, $A_F$ can never contains smooth
irreducible components to $\mathbb{C}$.
\end{theorem}

This fact was presented in \cite{Chau99} and can be reduced from
\cite{Cassou} (see also \cite{Chau2004} and \cite{Chau2007a} for
other refine versions). It can be used  in assuming that there is
a plane algebraic curve $E$,  composed of some irreducible
components parameterized polynomial maps $(\alpha(t),\beta(t))$
satisfying (1), such that the restriction
$$F:\mathbb{C}^2\setminus F^{-1}(E)\longrightarrow
\mathbb{C}^2\setminus E\eqno(2)$$ gives a unbranched covering.
Each of irreducible components of such a curve $E$ is a singular
curve approaching to $(\infty,\infty)\in \mathbb{P}^1\times
\mathbb{P}^1\supset \mathbb{C}^2$. In working with generic fibres
of $P$ we may use the following convenience, which was presented
in \cite{Chau2003} in a little different statement.

\begin{proposition}[See Theorem 1 in \cite{Chau2003}]\label{Singular} Suppose $F=(P,Q)$ is a Keller
map and $E$ is a plane algebraic curve in (2). If the line
$L_c:=\{ (c,t): t\in \mathbb{C}\}$ intersects transversally each
irreducible component of the non-proper set $E$, then $P=c$ is a
generic fiber of $P$ .
\end{proposition}

\vskip 0.3cm\noindent{\bf 3.}  Below we shall reduce from the
results introduced in the previous section some advance
estimations on the polar branches and the genus of generic fibres
of $P$ and $Q$

For each $c\in \mathbb{C}$  let us denote by $Pol(P,c)$ the union
of polar branches in the fiber $P=c$. The germ curve $Pol(P,c)$
can be realized as the inverse image $F^{-1}(V(c,R))$ for enough
large number $R>0$, where $V(c,R):=\{(c,t): |t|>R\}$.

\begin{theorem}\label{Pol} Suppose $F=(P,Q)$ is a Keller map. Then, the family of germ curves
$Pol(P,c)$, $c\in \mathbb{C}$, is equianalytical. Namely, for each
$c_0\in \mathbb{C}$ there exists a small dick $\Delta\subset
\mathbb{C}$ centered at $c_0$ and a number $R>0$ such that the
restriction
$$P:\bigcup_{c\in \Delta}F^{-1}(V(c,R))\longrightarrow \Delta$$
is an analytic fibration. In particular, the number of irreducible
components in $Pol(P,c)$ does not depended on $c\in \mathbb{C}$.
\end{theorem}

This fact was contained implicitly in Theorem 3.4 in \cite{Chau99}
, stated in terms of Newton-Puiseux expansions. For convenience,
we present here a proof by applying Theorem \ref{AF}.

\begin{proof}[Proof of Theorem \ref{Pol}] In view point of Theorem \ref{AF} we can assume that
that there is a plane algebraic curve $E$,  composed of some
irreducible components parameterized polynomial maps
$(\alpha(t),\beta(t))$ with $${\deg \alpha\over\deg \beta}={\deg
P\over\deg Q},$$ such that the restriction
$$F:\mathbb{C}^2\setminus F^{-1}(E)\longrightarrow \mathbb{C}^2\setminus
E$$ gives a unbranched covering.

Let $c_0\in \mathbb{C}$ be given. Since the components of $E$
approach to the point $(\infty,\infty)\in \mathbb{P}^1\times
\mathbb{P}^1\supset \mathbb{C}^2$, there exists a small dick
$\Delta\subset \mathbb{C}$ centered at $c_0$ and a number $R>0$
such that
$$(\Delta\times \{t\in \mathbb{C}: |t|>R\}) \cap E=\emptyset.$$ This
ensures that the restriction
$$F:F^{-1}(\Delta\times \{t\in \mathbb{C}: |t|>R\})\longrightarrow \Delta\times \{t\in \mathbb{C}: |t|>R\}$$
is a unbranched analytic covering. Note that $$F^{-1}(\Delta\times
\{t\in \mathbb{C}: |t|>R\})=\cup_{c\in \Delta}F^{-1}(V(c,R)).$$
This follows that $$P:\bigcup_{c\in
\Delta}F^{-1}(V(c,R))\longrightarrow \Delta$$ is an analytic
fibration.
\end{proof}

Let us denote by $g_h$ the genus of the generic fiber of a
primitive polynomial $h\in \mathbb{C}[x,y]$.

\begin{lemma}\label{genus} Suppose $F=(P,Q)$ is a Keller map.
If  $ \deg P\leq \deg Q$, then $ g_P \leq g_Q$.
\end{lemma}
\begin{proof}
If $F$ is invertible, of course  $g_P=g_Q=0$. Consider the
situation when $F$ is not invertible. We need  prove only that if
$\deg P\leq \deg Q$, then for each $t$ enough small a generic
fiber of $P$ can be topologically embedded into a generic fiber of
$P+tQ$. This ensures that $g_P\leq g_{P+tQ}$ for each $t$ enough
small that implies the desired conclusion.

Let $A_F$ be the non-proper value set of $F$. Note that $A_F$ is a
curve in $\mathbb{C}^2$ and the restriction $F:
\mathbb{C}^2\setminus F^{-1}(A_F)\longrightarrow
\mathbb{C}^2\setminus A_F$ gives a unramified covering. Replacing
$F$ by $F+p$ for a generic point $p\in \mathbb{C}^2$ if necessary,
we can assume that $(0,0)\in \mathbb{C}^2\setminus A_F$ and for
$\vert t\vert <\epsilon$ the lines $L_t$ given by $u+tv=0$
intersects transversally $A_F$. By Proposition \ref{Singular} the
late ensures that for $\vert t\vert <\epsilon$  the curve $P+tQ=0$
is a generic fiber of $P+tQ$.

Now, we will construct topological embeddings
$(P=0)\hookrightarrow (P+tQ=0)$ for $\vert t\vert <\epsilon$. To
do it, we can choose a box $ B:=\{ \vert u \vert < r; \vert u\vert
< s\}$ such that $(L_0\cap A_F) \subset B$ and $A_F\cap B$ is a
smooth manifold. Since the lines $L_t$ intersects transversal
$A_F$, by standard arguments we can modify the motion
$\phi_t(0,v):=(-tv,v)$ such that $\phi_t(A_F\cap L_0)\subset
A_F\cap L_t$  and $\phi_t: L_0\cap B\longrightarrow \phi_t(L_0\cap
B)\subset L_t\cap B$ are diffeomorphisms. Let $\Phi_t$ be the
lifting map $\Phi_t$ induced by the covering $F:
\mathbb{C}^2\setminus F^{-1}(A_F)\longrightarrow
\mathbb{C}^2\setminus A_F$. Then, $\Phi_t:F^{-1}(L_0\cap
B)\longrightarrow \{ P+tQ=0\}$ is an embedding of $F^{-1}(L_0\cap
B)$ into the fiber $P+tQ=0$. Since $(L_0\cap A_F) \subset B$, it
is easy to see that the fiber $P=0$ can be deformed diffeomorphic
to its subset $F^{-1}(L_0\cap B)$. So, we get the desire
embeddings.
\end{proof}

\vskip0.3cm\noindent{\bf 4.} Now, we are ready to prove the main
result.

\begin{proof}[Proof of Theorem 2] Let
$F=(P,Q)$ be a given Keller map. Assume that $P$ and $Q$ are
rational polynomials. We can assume in addition that the following
conditions holds:
\begin{enumerate}
\item[a)]  $\deg P <\deg Q$;

\item[b)] The curve $P=0$ is irreducible;

\item[c)] For generic $\lambda\in \mathbb{C}$ the curve $\lambda
P+Q$ is generic fiber of $\lambda P+Q$.
\end{enumerate}
Indeed,  if $\deg P=\deg Q$,  by the Jacobian condition we have
$P_+=cQ_+$ for a number $c\neq 0$, where $P_+$ and $Q_+$ are
leading homogenous components of $P$ and $Q$, respectively. Then,
$\deg P-cQ<\deg Q$ and $P-cQ$ is also rational by Lemma
\ref{genus}. So, we can replace $P$ by $P-cQ$. Further, in view of
Theorem \ref{AF} and Proposition \ref{Singular}, we can choose a
generic point $(a,b)\in \mathbb{C}^2$ such that the curve $P=a$ is
a generic fiber of $P$ and the curve $\lambda P+Q=\lambda a +b$ is
those of $\lambda P+Q$ for generic $\lambda\in \mathbb{C}$. So, we
can replace $F$ by $F-(a,b)$.

We will show that the fiber $P=0$ has only one polar branch. Then,
by Theorem 5 every fiber of $P$ also has only one polar branch.
This ensures that all fibres of $P$ are irreducible. Therefore, by
Razar's theorem  $F$ is invertible.

To do it, we regard the plane $\mathbb{C}^2$ as a subset of the
projective plane $\mathbb{P}^2$ and associate to $F=(P,Q)$ the
rational map $G_F:\mathbb{P}^2\longrightarrow \mathbb{P}^1$
defined by $G_F(x,y)=Q(x,y)/P(x,y)$ on the part $\mathbb{C}^2$.
The indeterminacy point set of $G_F$ consists of the set
$B:=F^{-1}(0,0)$ and a finite number of points in the line at
infinity $L_\infty$ of the chart $\mathbb{C}^2$. Note that the
restriction of $G_F$ to $L_\infty$ is equal to the zero, since
$\deg P <\deg Q$. By blowing-up we can remove indeterminacy points
of $G_F$ and obtain a blowing-up map $\pi:X\longrightarrow
\mathbb{P}^2$ and a regular extension $g_F: X\longrightarrow
\mathbb{P}^1$ over a compactification $X$ of
$\mathbb{C}^2\setminus B$.

Now, suppose $g_F:X\longrightarrow \mathbb{P}^1$ is such a regular
extension of $G_F$. Let $\mathcal{D}_\infty:=\pi^{-1}(L_\infty)$
and $\mathcal{D}_b:=\pi^{-1}(b)$, $b\in B$. $\mathcal{D}_\infty$
is a connected rational curve with simple normal crossing and its
dual graph is a tree. Each $\mathcal{D}_b$ is a copy of
$\mathbb{P}^1$. The divisor
$\mathcal{D}:=X\setminus(\mathbb{C}^2\setminus B)$ then is a
distinct union of $\mathcal{D}_\infty$ and $\mathcal{D}_b$.
Observe, for generic $\lambda\in \mathbb{P}^1$ the fiber
$g_F=\lambda$ is the closure in $X$ of the portion $\{(x,y)\in
\mathbb{C}^2\setminus B: \lambda P(x,y)+Q((x,y)=0\}$ of the fiber
$ \lambda P+Q=0$. Since $\deg P <\deg Q$ by Condition (a), in view
of Lemma \ref{genus} the polynomials $\lambda P+Q$ are rational.
Therefore, by Condition (c) generic fibres of $g_F$ are
irreducible rational curves. This means that $g_F:X\longrightarrow
\mathbb{P}^1$ is a $\mathbb{P}^1$-fibration over $\mathbb{P}^1$.

Now, we consider the fiber of $g_F$ over $\infty$, denoted by
$C_\infty$. Let us denote $D_\infty:=C_\infty\cap
\mathcal{D}_\infty$ and by $\Gamma$ the closure in $X$ of the
portion $\{(x,y)\in \mathbb{C}^2\setminus B: P(x,y)=0\}$. By the
conditions (a) and (b) $\Gamma$ is an irreducible rational curve
and $D_\infty$ contains at least the proper transform of the line
at infinity of $\mathbb{C}^2$. Furthermore, by the Jacobian
condition the multiplicity of $g_F$ on $\Gamma$ is equal to $1$.
Obviously, $C_\infty=D_\infty\cup \Gamma$ and $\Gamma$ intersects
$D_\infty$ at polar branches of the fiber $P=0$. By the well-know
fact (see for example \cite{Griff} and \cite{Ven}), that any
reducible fiber of a $\mathbb{P}^1$-fibration over $\mathbb{P}^1$
the can be contracted by blowing down to any its component of
multiplicity $1$, the fiber $C_\infty$ can be contracted to
$\Gamma$. In particular, $D_\infty$ is a blowing-up version of one
point and $C_\infty$ is a blowing-up version of $\mathbb{P}^1$.
Hence, the dual graphs of $C_\infty$ and $D_\infty$ are tree.
Therefore, the curve $\Gamma$ intersects transversally $D_\infty$
at a unique smooth point of $D_\infty$. This follows that the
fiber $P=0$ has only one polar branch.
\end{proof}

\vskip 0.3cm\noindent{\bf 5.} To conclude we would like to note
that it remains open the question {\it whether a Keller map
$F=(P,Q)$ with rational component $P$ is invertible}. In view of
Theorem 2 and its proof, if such a Keller map is not invertible,
then $\deg P <\deg Q$. In fact, it is possible to show that in
such a Keller map $\deg P$ does not divide  $\deg Q$ and the
fibres of $P$ have exactly two polar branches. We will return to
discuss on the question in a further paper.

 \vskip 0.3cm{\it Acknowledgements.} The first
versions of this paper contain some serious mistakes that was
pointed out by Pierete Cassou-Nogues. We would like to express our
thank to her for many valuable discussions. We also thank very
much Shreeram Abhyankar, Hyman Bass, Arno Van den Essen, Carlos
Gutierrez, L\^e D\~ung Tr\'ang, Walter Neumann, Mutsui Oka, Stepan
Orevkov and Nguyen Tien Zung for all valuable encouragements they
spend to us.

\end{document}